\documentclass{amsart}

\theoremstyle{definition}

\theoremstyle{remark}

\numberwithin{equation}{section}

\usepackage{amssymb,amsmath}

\def\o{\omega}
\def\u{\mu}
\def\us{\mu^*}

\def\p{\rho}

\def\xss#1#2{x_{#1}^{#2}(\u)}
\def\xsss#1#2{x_{#1}^{#2}(\us)}

\def\xk{x^k(\mu)}
\def\xk1{x^{k-1}(\mu)}

\def\du#1{\dfrac{d}{d#1}}

\def\cal#1{\mathcal{#1}}
\def\c{\circ}

\def\e#1{\exists{#1}}

\def\dfrac#1#2{\displaystyle{\frac{#1}{#2}}}

\def\raw{\Rightarrow}

\begin{document}

\title[Monotonicity of Finite Limit Cycles]{When will a One Parameter Family of Unimodal Maps Produce Finite Limit Cycles Monotonically with the Parameter?}

\author{John Taylor}
\address{Department of Mathematical Sciences, United Arab Emirates University, Al Ain, UAE }
\email{john.taylor@uaeu.ac.ae}
\thanks{}

\subjclass[2000]{Primary 37E05, 54H20; Secondary 37B40}

\date{September 24, 2007}

\keywords{One Parameter Family, Unimodal Map, Monotonicity, Kneading Sequence, Topological Entropy}

\begin{abstract}
In this note we consider a collection $\cal{C}$ of one parameter families of unimodal maps of $[0,1].$ Each family in the collection has the form $\{\mu f\}$ where $\mu\in [0,1].$  Denoting the kneading sequence of $\mu f$ by $K(\mu f)$, we will prove that for each member of $\cal{C}$, the map $\mu\mapsto K(\mu f)$ is monotone.  It then follows that for each member of $\cal{C}$ the map $\mu\mapsto h(\mu f)$ is monotone, where $h(\u f)$ is the topological entropy of $\mu f.$  For interest, $\mu f(x)=4\mu x(1-x)$ and $\mu f(x)=\mu\sin(\pi x)$ are shown to belong to $\cal{C}.$  This extends the work of Masato Tsujii [1].  
\end{abstract}

\maketitle

\section*{introduction}


Metropolis, Stein and Stein were among the first, to my knowledge, to study what are now called finite kneading sequences.  These were associated with one parameter families of interval maps which included $\mu f(x)=4\mu x(1-x)$ and $\mu f(x)=\mu\sin(\pi x)$ [2].  Computer studies strongly suggested a universal topological dynamics for a large class of such families, and many workers were quickly drawn to this emerging field of study. In the 1980's and 1990's there was intense interest in the behavior of the logistic map (affinely modified) in the setting of one complex dimension.  A central question was (essentially) under what circumstances would finite kneading sequences be monotone with the parameter, as this question was associated with the structure of the boundary of the Mandelbrot set.  This question was successfully addressed (in the special case of a real quadratic) in [1,5].  Here we address the question generally for a large class that includes the logistic map, and give a sufficient condition for the solution.  

All known proofs of this type apply to the case of a quadratic polynomial only and use complex analytic methods (holomorhic techniques) or depend on complex analysis (Compare [1], [5], [6], [7] [8]).  These methods are not used here.
       
Let $I=[0,1].$  Consider the collection $\cal{C}$ one parameter families of unimodal maps $\{\mu f\}$ where $x, \mu\in I,$ and $\mu f:I\rightarrow I$ with $\mu f$ at least $C^{(3)}$ in both $\mu$ and $x$.  Notice that since $\u f(x)$ is linear in $\u, $ it is $C^{\infty}$ in $\u.$  Denote the single critical point $c\in(0,1)$ and scale the map so that $f(c)=1,$ requiring that $f(0)=0=f(1).$  Then $\mu f(c)=\mu.$  Denote the $n^{th}$ iterate of $\mu f$ by $f_\mu^n(x)=(\mu f)\circ\cdots\circ(\mu f)\circ(\mu f)(x),$ where the composition is $n$-fold.
\bigskip

For any $x\in I$, the {\it orbit} of $x$ is the set $O(x)=\{f_\mu^n(x)| n\ge 0\}.$  We associate with $O(x)$ the word   
$\o(x) = \o_0\o_1\o_2\cdots$ with $\o_k\in\{L, C, R\}$ where words are formed as follows:

$$\o_k = \left\{ \aligned L, \quad \text{for} \quad f_\mu^k(x) < c\cr
C, \quad \text{for} \quad f_\mu^k(x) = c\cr
R, \quad \text{for} \quad f_\mu^k(x) > c. \endaligned \right.$$
$\o(x)$ is called the {\em itinerary} of $x$ under $f$.  Thus, $\o(x)$ will be infinite if and only if $O(x)$ is aperiodic.  We are interested in finite words.

In particular, we will be interested in studying the itinerary associated with $O(\u)$. This special itinerary is called the {\it kneading sequence} of $\mu f$, symbolized $K(\u f)=\o(\u)$.    
\bigskip

The following preliminaries are necessary for the statement that the map $\mu\mapsto K(\mu f)$ is monotone be meaningful;  that is, we need a total order on the kneading sequences.
\bigskip

It is possible to construct a total order $\prec$ on the set of all kneading
sequences, and more generally, on the set of all words made from the
alphabet $\{ L,C,R\},$ in such a way that it reflects the
order of the real line, in the sense that
 $x < y$ implies $\o(x) \preccurlyeq \o(y)$, is defined as follows:  First define $L\prec C\prec R.$  Then, if
$A=\{a_k\}\ne B=\{b_k\},$ let $N$ be the smallest index for which $a_N\ne b_N$
and let $\p_{N-1}$ be the number of $R's$ in the word $a_1\cdots a_{N-1}.$
Then define $A\prec B$ if $a_N\prec b_N$ and $\p_{N-1}$ is even or if $a_N\succ b_N$ and $\p_{N-1}$ is odd.

This order it sometimes referred to as the {\em parity-lexicographic} order.  
The intuition derives from the fact that for $x\in[c,1],$ $\u f$ is orientation reversing, that is, $x<y$ implies that $f(x)>f(y).$  In order that the ordering on the words be consistent with the order in the real numbers, this reversal of orientation is accounted for in the manner just described.

A word $\o$ is called {\it maximal} (or {\it shift-maximal}) provided it
is greater (in the parity lexicographic order) than all of its shifts,
where, as usual, the shift operator $\sigma$ is defined by the
action $\sigma(\o) = \o_1\o_2 \o_3\cdots$ on the word $\o=\o_0\o_1\o_2\o_3\cdots .$  Shift maximal words correspond to periodic orbits. 

In kneading theory there are several versions of an ``intermediate value
theorem''.  This type of theorem is fundamental in that it relates abstract words to the behavior of dynamical systems.  That is, it connects the set of kneading sequences ordered by the relation $\prec$ and the parameter space (an interval in the real line)
with the usual order.  The following version is essentially that found in
[3]: 
\bigskip

{\bf Theorem A}\ Let $\{\mu f \}$ be any one parameter family of
$C^1$ unimodal maps.  If $\mu_1 < \mu_2$ are two parameter values with
corresponding kneading sequences $K(\mu_1 f)\prec K(\mu_2 f)$, and if
$\o$ is any shift-maximal sequence with the property that
$$K(\mu_1 f)\prec \o\prec K(\mu_2 f),$$
then there exists a $\mu$ such that $\mu_1 < \mu < \mu_2$ and $\o =
K(\mu f)$.
\bigskip

Since $\mu f$ double covers $[0, \u]$ in such a way that $\mu f([0,c])=[0, \u]=\mu f([c,1]),$ the functions $(\mu f)_{\o}^{-1},\ \o\in\{L, R\}\ \text{ have the action}\ (\mu f)_{L}^{-1}([0, \u])=[0,c]\ \text{and}\ (\mu f)_{R}^{-1}([0, \u])=[c,1].$

For all $n\ge 0,$ let $G_n(\mu)$ denote the graph of $f_\mu^n.$

By the chain rule, $\dfrac{d}{dx}f^n_\mu(x)=\prod_{k=0}^{n-1}\mu f'[f_\mu^k(x)].$
Therefore,  if $x$ is an extreme point of $G_n(\mu)$ then there exists $k, 0\le k\le n$ such that $f^k(x)=c$.
\bigskip

{\bf Definition}\ 
For $\mu$ fixed, define $\xss{}{0}=c,$ and for $1\le k\le n-1,$ denote by the symbol $x^k_\o$ any $k^{th}$ preimage of $c,$ specifically
$$
f^{-k}_\mu(c)=\{\xss{\o}{k}:=(\mu f)_{\o_1}^{-1}\circ \cdots\circ (\mu f)_{\o_{k}}^{-1}(c)\ |\ P=\o_1 \o_2 \dots \o_k\quad \o_i\in\{L,R\}\}.
$$

Fixed points of the  functions $\xss{\o}{k}$, which, given $f$, are functions of $\mu$ alone, will be central in what follows 

Denote the graph $\Gamma f^k_\u$ by $G_k(\mu)$.



For any $C^3$ function $\psi$, let $S(\psi)= \dfrac{\psi'''}{\psi'}-\dfrac{3}{2}\left(\dfrac{\psi''}{\psi'}\right)^2.$ 
A simple computation reveals that if $\phi$ is also a $C^3$ function with $S(\phi)>0$ and $S(\psi)>0$, then $S(\psi\circ \phi)>0.$

\section*{main section}

Here we prove that for each member of $\cal{C}$, the map $\mu\mapsto K(\mu f)$ is monotone.
\bigskip

All families in $\cal{C}$ have the following properties:
\bigskip

1)\ For each $\mu$ there exists a unique fixed point for $\mu f$ in
$(0,1)$.  
\bigskip

2)\ For each fixed $\mu$ and for all $n\ge 1,$ $f^n_\mu$ has at most one attracting periodic orbit,
and $O(\u)$ is asymptotic to this attracting periodic orbit.
\bigskip

3)\  $S[(\mu f)_{\o}^{-1}]>0$  for all $\mu,$ where $\o\in\{L, R\}.$
\vskip 20pt

{\bf Remarks}
\bigskip

(i)\ Concave maps, for example, have property 1.
\bigskip

(ii)\ It is known that if $S(f)<0$ for all $x$, then property 2 holds. [4]
\bigskip

{\bf Lemma} Assume that $\o=\o_1\o_2\cdots \o_{k-1}=K(\u f)$ and that $f_\u^k(c)=c$ has primitive period $k,$ so that by continuity there is an open set of parameter values for which the composition
$$(\mu f)_{\o_1}^{-1}\circ \cdots\circ (\mu f)_{\o_{k-1}}^{-1}(c)$$
 is defined.  If $\u$ is such that $f_{\u}^k(c)=c,$  then
$$ x^{k-1}_{\o}(\mu)=\mu\ \Rightarrow\ f_{\u}^k(c)=c.$$
\smallskip

{\bf Proof}\quad  $$x^{k-1}_{\o}({\mu})=({\mu} f)_{\o_{1}}^{-1}\c\cdots\c ({\mu} f)_{\o_{k-1}}^{-1}(c)={\mu}\ \Rightarrow$$
$$c=f^{k-1}_{\mu}\left[
x^{k-1}_\o({\mu})
\right]
=f^{k-1}_{\mu}\left[
({\mu} f)_{\o_{1}}^{-1}\c\cdots\c ({\mu} f)_{\o_{k-1}}^{-1}(c)
\right]
=f^{k-1}_{\mu}({\mu})=f^{k-1}_{\mu}\left[
{\mu} f(c)
\right]=f_{\mu}^k(c).$$
$\quad\square$

%
\bigskip

{\bf Remarks}
\bigskip

(i)\ A super stable point of primitive period $n$ occurs in association with the equation $\xss{\o}{n-1}=\mu$,  where $K(\mu f)=\o$.
\bigskip

(ii)\ The trajectories of distinct preimages can never intersect.
\bigskip

It follows from the implicit function theorem that, for all $n\ge 1,\ 0\le k\le n-1,$ level functions of order $k$ exist so long as the intersection of $G_{n}(\mu)$ with the line $y=c$ exists.

But this intersection exists for all $\mu>\mu^*$, where $\mu^*$ is the parameter value with the property that, for $1\le k\le n,$ 
$x_{\o}^{k-1}(\mu^*)=\mu^*$;  for then, $f^{k}_{\mu^*}(c)=c $ in $G_{n}(\mu)$ by the lemma.  Therefore, so long as $\mu^*$ is unique with the above property, we see that for all $\mu>\mu^*$, the intersection of $G_{n}(\mu)$ and the line $y=c$ persists, and so, the level functions $\xss{\o}{n-1}$ exist on a connected domain.
Further, a certain number of these $\xss{\o}{n-1}$ will have fixed points.  The number is known to be  $$\displaystyle{\dfrac{1}{2n}\sum\mu(d)2^{n/d}}$$
where the sum is taken over all odd square free divisors of $n.$ [5]
\bigskip

{\bf Theorem}\quad For each member of $\cal{C}$, the map $\mu\mapsto K(\mu f)$ is monotone.
\bigskip

{\bf Proof}\quad The proof is by strong induction.  First, notice that $\xss{R}{1}$ exists on the connected domain $[c,1]$.  Since we assume that $S(\xss{R}{1})>0$, $\du{\mu}\xss{R}{1}$ cannot have a positive local maximum.  Therefore,
$$
\e{!\us}\left[\xsss{R}{1}=\us\right]\ \Rightarrow\  \e!\us\left[f_{\us}^2(c)=c\right]\quad\text{for}\quad K(\us f)=RC\quad \text{by the lemma}.
$$
In other words, $\e!\mu_{\o}\ [f^1_{\mu_{\o}}(c)=c]$ with $\o=K(\mu_\o f) \raw\  \e!\mu_{\tau}\ [f^2_{\mu_{\tau}}(c)=c]$ with $\tau=K(\mu_{\tau} f).$  Here $\o=C$ and $\tau=RC$.
\smallskip

Assume that for $1\le k\le n,$ and for all $\o=K(\mu f)$ (with the length of $\o$ not exceeding $n$),  $\e{!\mu_\o}\ [x_{\o}^{k-1}(\mu_\o)=\mu_\o]$, that is $\e!\mu_\o\ [f_{\mu_\o}^k(c)=c]$ with $\o=K(\mu f).$
\smallskip

Since $\e{!\mu_\o}\ [x_{\o}^{k-1}(\mu_\o)=\mu_\o]$, that is, $\e!\mu_\o\ [f_{\mu_\o}^k(c)=c],$ dom$(\xss{\o}{k})$ is connected.
\smallskip

If  $\xss{\o}{k}$ has a fixed point, that is, if $P=K(\mu f)$ for some $\mu$, then $f_{\mu}^{k+1}(c)=c$ when $\xss{\o}{k}=\mu.$
\smallskip

But $S[\xss{\o}{k}]>0\ \raw\ \e{!\mu}\ [x_{\o}^{k}(\mu)=\mu],$ that is, $\e!\mu\ [f_{\mu}^{k+1}(c)=c]$ with $\o=K(\mu f).$
\smallskip

In particular, $\e!\mu_\o\ [f_{\mu_\o}^n(c)=c]$ with $\o=K(\mu_{\o} f) \raw\  \e!\mu_{\tau}\ [f_{\mu_{\tau}}^{n+1}(c)=c]$ with $\tau=K(\mu_{\tau} f).$\quad$\square$
 \bigskip

{\bf Remarks}\  (1) One computes that $S[(\mu f)^{-1}_\o]>0,\ \o\in\{L,R\}$ when $\mu f(x)=4\mu x(1-x)$ and $\mu f(x)=\mu\sin(\pi x).$
\bigskip 

\hskip 0.675in (2) The topological entropy of maps in the class $\cal{C}$ is evidently monotone with the parameter.  This is because orbit production for these never decreases.
\bigskip

\end{document}